\documentclass{amsart}
\usepackage{a4}
\usepackage{amssymb}

\newtheorem{Th}{Theorem}[section] 
\newtheorem{Lem}[Th]{Lemma} 

\newtheorem{Rem}[Th]{Remark}

\newtheorem{claim-num}{Claim}

\numberwithin{equation}{section}
\renewcommand{\theequation}{\thesection.\arabic{equation}}

\def\gl#1{\operatorname{GL}(#1)}

\def\aut#1{\operatorname{Aut}(#1)}

\def\wid{\operatorname{wid}}
\def\codim{\operatorname{codim}}
\def\id{\operatorname{id}}

\def\inv{^{-1}}
\def\str#1{\langle#1\rangle}

\def\f{\varphi}
\def\a{\alpha}
\def\b{\beta}
\def\s{\sigma}
\def\vk{\varkappa}

\def\Q{\bold Q}

\def\cB{{\mathcal B}}

\def\N{\mathbf N}

\def\Q{\mathbf Q}
\def\R{\mathbf R}

\renewcommand{\le}{\leqslant}

\begin{document}

\begin{abstract}
Recently George Bergman proved that the symmetric
group of an infinite set possesses the following
property which we call by the {\it universality of
finite width}:  given any generating set $X$ of the
symmetric group of an infinite set $\Omega,$ there is
a uniform bound $k \in \N$ such that any permutation
$\sigma \in \text{Sym}(\Omega)$ is a product of at
most $k$ elements of $X \cup X^{-1},$ or, in other
words, $\text{Sym}(\Omega)=(X^{\pm 1})^k.$ Bergman also
formulated a sort of general conjecture stating that
`the automorphism groups of structures that can be put
together out of many isomorphic copies of themselves'
might be groups of universally finite width and
particularly mentioned, in this respect,
infinite-dimensional linear groups. In this note we
confirm Bergman's conjecture for infinite-dimensional
linear groups over division rings.  \end{abstract}

\title[Infinite-dimensional general linear groups...]{Infinite-dimensional general linear
groups are groups of universally
finite width}
\author{Vladimir Tolstykh}
\address{Department of Mathematics\\ Yeditepe University\\
34755 Kay\i\c sda\u g\i \\
Istanbul\\
Turkey}
\email{vtolstykh@yeditepe.edu.tr}
\date{March 13, 2004}
\maketitle

Let $G$ be a group and $S$ a generating set of $G.$
Recall that given an element $x$ of $G,$ we call by the {\it length}
$|x|_S$ of $x$ relative to $S$ the minimal
natural number $k$ such that $x$ is a product of $k$
elements from $S \cup S\inv.$ Then the {\it width} $\wid(G,S)$ of $G$
relative to $S$ is the supremum of all numbers
$|x|_S,$ where $x$ runs over $G.$ In particular,
$\wid(G,S)$ is either a natural number, or $\infty.$

It might be said that in most of cases one expects to
find in a given group both generating sets that
provide finite widths and generating sets that provide
infinite widths.  For instance, if an infinite group
$G$ is generated by some finite set $S,$ then
$\wid(G,S)$ is infinite, while, on the other hand,
$\wid(G,G)=1.$

However, following Bergman \cite{Berg},
we can formally set an axiom
saying that
\begin{quote}
for every generating set $S$ of a
group $G$ the width of $G$ with
respect to $S$ is finite
\end{quote}
thereby defining the class of groups of, say, {\it
universally finite width.} Most obvious examples of
groups satisfying this axiom are finite groups. Known
examples of infinite groups of universally
finite width are all rather non-trivial (and
let us note in passing that all
these groups are {\it uncountable}). First results in
this direction were obtained by Shelah in 1980 in his
paper \cite{Sh}, where he constructed such an
uncountable group $G$ that $G$ had width at most 240
relative to any generating set.

A similar result obtained last year by Bergman and stating
that the symmetric group of an infinite set is a group
of universally finite width had a considerable impact.
Without exaggeration, Bergman's result charmed many a
mathematician.  Soon afterwards, it has been
understood that the automorphism groups of doubly
homogeneous ordered sets (Droste-Holland, \cite{DrHo}),
the automorphism group of $\R$ as
a Borel space and the homeomorphism groups of certain
topological spaces (such as, for instance, the Cantor
discontinuum, $\Q$ and the set of irrational numbers; Droste-G\"obel, \cite{DrGo}) are groups of universally
finite width. These results confirmed a general
conjecture expressed by Bergman in the first
version of his preprint \cite{Berg} as follows:
\begin{quote}
``It seems likely that analogues of our result should be
provable for the automorphism groups of other
structures that can be put together out of many
isomorphic copies of themselves; e.g., the
automorphism groups of $\R$ as a Borel space, of a free
group on infinitely many generators, or of some of the
structures whose automorphism groups are studied in
[*], ..., [*].''
\end{quote}
Among the automorphism groups `studied in [*], ...,
[*]' are, in particular, infinite-dimensional general
linear groups. The aim of this note is to justify for
the groups in question the inspring conjecture of Bergman's. We
also suggest a way of generalization of our result to
the automorphism groups of free algebras of infinite
ranks.

Bergman's preprint \cite{Berg} is founded on the ideas
from the famous paper by Macpherson-Neumann \cite{MN}
on subgroups of infinite symmetric groups and very
skillfully develops these ideas -- particularly, those
ideas from \cite{MN} that are aimed at determination of the confinalities
of symmetric groups.  Some general results on confinalities of
infinite-dimensional general linear groups were
obtained by Macpherson in the paper \cite{MacPh}, as a
natural continuation of the study of infinite
symmetric groups.  It turned out that Macpherson's
ideas can be also adapted to dealing with the widths
of infinite-dimensional general
linear groups.

The author would like to express his sincere gratitude to
Valery Bardakov for helpful and fruitful conversations
on the subject of the present note.

\section{$\gl V$ as a power of some conjugacy class}

\def\WidofCpi{28}

By $V$ we shall denote a left infinite-dimensional vector
space over a division ring, the dimension
of $V$ will be denoted by $\vk.$ $G$ stands
for the general linear group $\gl V$ of $V,$ the group of all invertible
linear transformations from $V$ into itself.  We shall
use standard permutation notation, extending it, like in
\cite{Berg}, to arbitrary subsets of $G.$ Thus, if
$Y$ is a subset of $G$ and $W$ is a subspace of $V,$
$Y_{(W)}$ is the set of all elements of $Y$ that fix $W$ pointwise
and $Y_{\{W\}}$ is the set of all elements of $Y$ that
fix $W$ setwise. Any notation like $Y_{*_1,*_2}$
means the set $Y_{*_1} \cap Y_{*_2}.$

Following \cite{MacPh}, we call a subspace
$U$ of $V$ by a {\it moietous} subspace if
$$
\dim U = \codim U =\vk =\dim V.
$$
The notion of a moietous subspace
stems from the notion of a {\it moiety}
of an infinite set $I,$ such a subset $J$ of $I$ that $|J|=|I \setminus J|.$

We fix throughout an involution $\pi^*$ of $V$ that
interchanges two moietous subspaces $U,U'$ of $V$ with
$$
V = U \oplus U',
$$
that is,
$$
\pi^* U = U'.
$$

Our main result stating that
$\gl V$ is a group of universally
finite width will be obtained
as a consequence of the following
statement.

\begin{Th} \label{Wid1-of-GL}
The width of $\gl V$ relative
to the conjugacy class $C(\pi^*)$ of $\pi^*$
is at most $\WidofCpi.$
\end{Th}

Thus,
$$
\gl V =C(\pi^*)^{\WidofCpi}
$$
which explains the title of the
section.

Theorem \theTh, in turn, is based on obtaining of an
estimate of the width of $\gl V$ relative to the
union of a pair of naturally defined subgroups,
introduced by Macpherson in \cite{MacPh}.

\def\MacLemWid{7}

\begin{Lem} \label{Mac}
Suppose that $U_1,U_2,W$ are moietous subspaces of $V$ such
that
$$
V = U_1 \oplus U_2 \oplus W.
$$
Let
$$
H_1 = G_{(U_2),\{U_1+W\}} \text{ and }
H_2 = G_{(U_1),\{U_2+W\}}.
$$
Then the width of $\gl V$ relative to
the set $H_1 \cup H_2$ is at most $\MacLemWid$:
$$
\wid(\gl V, H_1 \cup H_2) \le \MacLemWid.
$$
\end{Lem}

\begin{proof}
Our proof is based heavily on Macpherson's proof of
the fact that $\gl V$ is generated by $H_1 \cup H_2$
(see \cite{MacPh}, the proof of Proposition 2.2).

Let $\cB_1,\cB_2$ and $\cB_W$ be bases of
the subspaces $U_1,U_2$ and $W$ respectively.
Write $\cB$ for $\cB_1 \cup \cB_2 \cup \cB_W.$

Let $\cB'$ be an arbitrary ordered base of $V.$ Following
Macpherson, we are going to prove existence of $\le \MacLemWid$ elements
of $H_1 \cup H_2$ whose product, say, $\s$ takes $\cB'$ to
$\cB$ (being considered as an ordered,
too.) The hardest part in Macpherson's proof is to
find an element $\rho \in \str{H_1 \cup H_2}$ such
that
$$
\rho \cB' \supseteq \cB_1 \cup \cB_W.
$$
Once such an automorphism $\rho$ is found, it is
required to follow it by some element
from $H_2 H_1 H_2$
to move $\rho \cB'$ onto $\cB$ (see
\cite{MacPh} for the details.) Hence
$$
|\s| \le |\rho| + 3,
$$
where $|\f|$ denotes the length
of an element $\f$ of $\gl V$ relative
to $H_1 \cup H_2.$ Thus, it remains
to prove existence of $\rho$
whose length is at most $4$.

To do this, we are going to change slightly the
corresponding part of Macpherson's proof, since in its
course he several times takes elements from $\str{H_1
\cup H_2}$ (not from $H_1 \cup H_2$ as it would
suit us), and since, by understandable
reasons, he does not care as to uniform boundedness of
lengths of such elements.

Consider the natural projection $\pi$ of $V$
onto $U_1.$ Then the set $\pi \cB'$ generates
$U_1$ and hence an appropriate element $\gamma_1$ of
$H_1$ takes it onto a superset of $\cB_1$:
$\gamma_1 \pi \cB' \supseteq \cB_1.$
Note that $\gamma_1 \pi = \pi \gamma_1$ and
then
$$
\pi \gamma_1 \cB' \supseteq \cB_1.
$$
Suppose that
$$
\cB_1 = \{u_{i1} : i \in I\},
$$
where index set $I$ is of cardinality
$\vk.$ For each $i \in I$ we find
an element $x_i \in \gamma_1 \cB'$
whose projection is $u_{i1}$:
\begin{equation}
x_i = u_{i1}+v_i,
\end{equation}
where $v_i \in U_2 + W.$
Such a moiety $J$ of $I$ can be found
that the subspace $L=\str{v_i~:~i \in J}$
has codimension $\vk$ in $U_2+W.$
If so, we act on $L$ by an element
$\gamma_2 \in H_2$ which moves $L$ onto a subspace $W$
of codimension $\vk$ in $W$. This done,
we look for an element $\alpha$ of $H_1$ which moves
the subspace $\gamma_2 L$ into $U_1$ and
fixes pointwise the elements
of the part $\{u_{i1} : i \in J\}$
of the base $\cB_1$ and, furthermore,
such that the subspace
$$
\alpha \gamma_2 L +\str{u_{i1} : i \in J}
=\alpha \gamma_2 \str{u_{i1}, v_i : i \in J}
$$
is a moietous subspace of $U_1.$ By
(\theequation) we have that the
set
\begin{equation}
\{\alpha \gamma_2 x_i : i \in J\}
\end{equation}
also generates a moietous subspace of
$U_1.$ Let now $\cB_W^0$ be a moiety
of $\cB_W.$ It is possible to find
a $\beta \in H_1$ which maps
the set of vectors in (\theequation)
onto $\cB_1 \cup \cB_W^0$:
$$
\beta\{\alpha \gamma_2 x_i : i \in J\}=\cB_1 \cup \cB_W^0.
$$
Write $\gamma_3$ for the product $\beta \alpha$
from $H_1.$ Finally, we observe that
there is an automorphism $\gamma_4$ from $H_2$ which maps $\cB_W^0$
onto $\cB_W.$ Summing up, we see that
we can take as $\rho$ the automorphism
$$
\gamma_4 \gamma_3 \gamma_2 \gamma_1
$$
whose length is at most $4,$ as
required.
\end{proof}

\begin{proof}[Proof of Theorem \ref{Wid1-of-GL}]
Let $L_1,L_2,M$ be three moietous subspaces
such that $V$ is their direct sum:
$$
V=L_1 \oplus M \oplus L_2.
$$
Let $I$ be an index set of cardinality
$\vk$ and
$$
\{a_i :i \in I\},\quad \{a_i^* : i \in I\}
$$
be bases of $L_1$ and
$$
\{b_i : i \in I\}
$$
a base of $L_2.$ Involutions $\pi_1,\pi_2$
which both fix $M$ pointwise and act on the bases
$\{a_i\}$ and $\{a_i^*\}$ as follows
\begin{align*}
&\pi_1 a_i =b_i, \qquad \forall i \in I\\
&\pi_2 a_i^* =b_i
\end{align*}
are conjugates of $\pi^*.$ We have therefore
that
$$
\pi_2 \pi_1 a_i =a_i^*\quad \forall i \in I.
$$
Let now $\a$ denote the automorphism
of the vector space $L_1$ that takes
the base $\{a_i\}$ onto the base $\{a_i^*\}.$
Suppose that for all $i \in I$
$$
\a\inv a_i^* =\sum_j \b_{ij} a_j^*.
$$
We then have
$$
\pi_2 \pi_1 b_i =\pi_2 a_i =\pi_2(\a\inv a_i^*)
=\pi_2(\sum_j \b_{ij} a_j^*)=\sum_j \b_{ij} b_j
$$
for all $i \in I.$ We see that the action of $\pi_2\pi_1$
on $L_2=\str{b_i : i \in I}$ is isomorphic to the
action of $\a\inv$ on $L_1=\str{a_i^* : i \in I},$
or, informally, one can write that
$$
\pi_2 \pi_1 = \a \oplus \id\oplus \a\inv.
$$

Extending the principle
of the construction of $\pi_2\pi_1,$ one can
represent as a product of two conjugates of $\pi^*$ any
automorphism of $V$ of the form
\begin{equation}
\bigoplus_{n \in \N} \a \oplus \id \oplus
\bigoplus_{n \in \N} \a\inv,
\end{equation}
where the latter direct sum corresponds
to a direct sum of {\it moietous} subspaces
and $\a$ is the isomorphism type of an
automorphism of one of these subspaces.

Let us consider two automorphisms $\s_1,\s_2$ of $V$
of the form (\theequation):
\begin{alignat*}9
& \s_1 &= &\id &\oplus (&\a     &\oplus &\a\inv &\oplus &\a     &\oplus &\a\inv &\oplus &\ldots )     &\oplus  &\id,\\
&\s_2  &= &\a  &\oplus (&\a\inv &\oplus &\a     &\oplus &\a\inv &\oplus &\a     &\oplus &\ldots )     &\oplus  &\id
\end{alignat*}
(both constructed over the {\it same} decomposition
of $V$ into a countable infinite direct sum of moietous subspaces.)
Then
\begin{equation}
\s_1 \s_2 =\a \oplus \bigoplus_{n \in \N} \id \oplus \id
\end{equation}
is a product of four conjugates of $\pi^*$
(cf. with the proof of Claim 4.11 from
the paper \cite{Tho} by Thomas.)

Let $U_1,U_2,W$ be subspaces of $V$ with
$$
V = U_1 \oplus W \oplus U_2.
$$
The formula (\theequation) demonstrates
that each element of the subgroup
$G_{(U_1),\{U_2+W\}}$ (resp. the
subgroup $G_{(U_1),\{U_1+W\}}$)
is a product of at most four
conjugates of $\pi^*.$
Then by Lemma \ref{Mac} any
element of $\gl V$ is a product
of at most
$\MacLemWid \cdot 4 = \WidofCpi$ conjugates
of $\pi^*.$
\end{proof}

\begin{Rem} \em
An upper bound for the width of $\gl V$ relative to
the conjugacy class of $\pi^*$ we have obtained is of
course very inaccurate. This is caused by the use of
Macpherson's lemma, but this is not unavoidable. On
the other hand, namely Macpherson's lemma makes the
proof as easy as possible, and, moreover, we are going
to apply Macpherson's lemma once again in a more
serious circumstances.  \end{Rem}

\section{$\gl V$ is a group of universally finite width}

\begin{Th} \label{Wid-of-GL-Is-Finite}
Let $X$ be any generating
set of $\gl V.$ Then the
width of $\gl V$ with
respect to $X$ is finite.
\end{Th}

\begin{proof}
We start with a linear analogue
of the result that plays a crucial
role in Bergman's proof of the
fact that any symmetric group
is a group of universally finite
width.

\begin{Lem} \label{BergLem}
There exist a power $X^m$ of $X$ and a decomposition
$V=U \oplus W$ of $V$ into a sum of moietous subspaces
such that the set
$$
(X^m)_{\{U\},\{W\}}
$$
induces the group $\gl U$ on $U.$
\end{Lem}

\begin{proof} (By Bergman-Macpherson,
cf. \cite[Lemma 4]{Berg}, \cite[proof of Theorem 3.1]{MacPh}).
Let
$$
V = \bigoplus_{k \in N} L_k
$$
be decomposition of $V$ into a
sum of moietous subspaces.
Write
$$
L_k^* = \bigoplus_{i \ne k} L_i
$$
for all $k \in \N.$

If for some pair $(X^k,L_j)$ we have that
$$
(X^k)_{\{L_j\},\{L_j^*\}} \text{ induces $\gl{L_j}$ on $L_j$}
$$
then the conclusion of Lemma is true. Suppose
otherwise. Then, in particular, for all $k$
$$
(X^k)_{\{L_k\},\{L_k^*\}} \text{ does not induce $\gl{L_k}$ on $L_k$}.
$$
Then for each $k \in \N$ we can find $\s_k \in \gl{L_k}$ such
that
\begin{quote}
$\sigma_k$ does not equal to the
restriction on $L_k$ of any element from $(X^k)_{\{L_k\},\{L_k^*\}}$.
\end{quote}
Set
$$
\sigma = \bigoplus_{k \in \N} \sigma_k.
$$
Since $\gl V = \bigcup_k X^k,$ we have $\sigma \in X^j$
for some $j \in \N.$ It is clear, however, that
$$
\sigma \in (X^j)_{\{L_j\},\{L_j^*\}}.
$$
But then the restriction of $\sigma$ on $L_j$ is $\sigma_j,$ a
contradiction.
\end{proof}

Let $X^m,U,W$ satisfy the conclusion
of the Lemma. We are keeping in mind application
of Lemma \ref{Mac} and, therefore,
we are going first to wipe out any
traces of non-triviality of action
of $(X^m)_{\{U\},\{W\}}$ on $W.$

For this purpose we take an involution $\pi \in \gl V$
which is  conjugate to $\pi^*,$ fixes $U$ setwise
and fixes $W$ pointwise. Then by Theorem
\ref{Wid1-of-GL} the set of automorphisms
\begin{equation}
Y=\{ \s_1 \pi \s_1\inv \ldots \s_{\WidofCpi} \pi \s_{\WidofCpi}\inv :
\s_1,\ldots,\s_{\WidofCpi} \in (X^m)_{\{U\},\{W\}} \}
\end{equation}
is the group $G_{\{U\},(W)},$ since any automorphism
$\s_k \pi \s_k\inv$ acts trivially on $W$ due to
triviality of the action of $\pi$ on $W$ (cf. \cite[the proof
of Lemma 3.2]{MacPh}). Clearly, $Y$ is contained
in some power of $X.$ Now by Lemma \ref{Mac} for the conjugate
set $\rho Y \rho$ by a suitable involution $\rho$
we have that
$$
\wid(\gl V,Y \cup \rho Y \rho) \le \MacLemWid.
$$
Let $l$ be such a natural number that
$$
\rho Y \rho \subseteq X^l.
$$
Then, evidently,
$$
\gl V =X^{\MacLemWid l}.
$$
\end{proof}

\begin{Rem}
\em The reader may have a feeling that
the results obtained can be true
for a wider class of automorphism
groups. Let us outline one such possible
generalization (the interested reader may refer also
to the paper \cite{DrGo} by Drose-G\"obel, where the authors
suggest an axiomatic system(s) whose
models among infinite permutation
groups are groups of universally
finite width.)

{\bf Definition.} Let $\bf V$ be
a variety of algebras. We shall
call $\bf V$ by a {\it BMN-variety}
(Bergman-Macpherson-Neumann variety)
if for every free algebra $F$ of
infinite rank from $\bf V,$ given
any moietous subalgebras $U_1,U_2,W$ of $F$
with
$$
F = U_1 * U_2 * W
$$
where $*$ denotes the operation
of free product, we have that the automorphism
group $G=\aut F$ of $F$ is generated
by the subgroups
$$
G_{(U_2),\{U_1 * W\}} \text{ and }
G_{(U_1),\{U_2 * W\}}.
$$
Examples of BMN-varieties are therefore
the variety of all sets with no structure
(Macpherson-Neumann, \cite{MN}) and
any variety of vector spaces over
a fixed division ring (Macpherson, \cite{MacPh}).

Using the scheme of
the proof of Theorem \ref{Wid-of-GL-Is-Finite},
one can prove the following result.

\begin{Th}
Let $\bf V$ be a BMN-variety of
algebras. Then the automorphism group
$G=\aut F$ of a free algebra $F \in \bf V$
of infinite rank is a group
of universally finite width
if and only if the width
of $\aut F$ relative
to the set
$$
G_{(U_2),\{U_1 * W\}} \cup
G_{(U_1),\{U_2 * W\}}.
$$
over some moietous subalgebras
$U_1,U_2,W$ of $F$ such that
their free product is $F$ is finite.
\end{Th}

\begin{proof}
The necessity part is obvious. To prove
the sufficiency, we prove for $\aut F$
analogues of Theorem \ref{Wid1-of-GL}
and Theorem \ref{Wid-of-GL-Is-Finite}
(no serious changes required).
\end{proof}
\end{Rem}

\end{document}